\documentclass[12pt]{amsart}

 \usepackage{amsfonts,graphicx,amsmath,amsthm,amsfonts,amscd,amssymb,amsmath,latexsym,multicol,mathrsfs}
\usepackage{epsfig,url}
\usepackage{flafter}
\usepackage{fancyhdr}
\usepackage{hyperref}

\hypersetup{colorlinks=true, linkcolor=black}

\makeatletter

\def\jobis#1{FF\fi
  \def\predicate{#1}%
  \edef\predicate{\expandafter\strip@prefix\meaning\predicate}%
  \edef\job{\jobname}%
  \ifx\job\predicate
}

\makeatother

\if\jobis{proposal}%
\else
\fi

 \usepackage[matrix, arrow]{xy}

\DeclareMathOperator{\Supp}{Supp}

\DeclareMathOperator{\Spec}{Spec}

\DeclareMathOperator{\Div}{div}
\DeclareMathOperator{\Chow}{Chow}

 \numberwithin{equation}{subsection}
 \numberwithin{footnote}{subsection}

\newtheorem{thm}{Theorem}[section]

 \newtheorem{cor}[thm]{Corollary}
 \newtheorem{lem}[thm]{Lemma}
 \newtheorem{prop}[thm]{Proposition}
 \newtheorem{conj}[thm]{Conjecture}

{
    \newtheoremstyle{upright}%
        {8pt plus2pt minus4pt}%
        {8pt plus2pt minus4pt}%
        {\upshape}%
        {}%
        {\bfseries\scshape}%
        {}%
        {1em}%
        {}%
\theoremstyle{upright}

 \newtheorem{defn}[thm]{Definition}

}

 \newcommand{\PP}{\mathbb P}
 
 \newcommand{\Q}{\mathbb Q}
 \newcommand{\R}{\mathbb R}

 \newcommand{\OO}{\mathcal{O}}

 \newcommand{\bir}{\dashrightarrow}
 \newcommand{\lin}{\sim}
 \newcommand{\num}{\equiv}
 \newcommand{\isom}{\cong}
 \newcommand{\rddown}[1]{\left\lfloor{#1}\right\rfloor} 

\begin {document}

\title{F\lowercase{inite generation of the log canonical ring for 3-folds in char p}}
\author{J\lowercase{oe} W\lowercase{aldron}}

\begin{abstract}
We prove that the log canonical ring of a projective klt pair of dimension $3$ with $\Q$-boundary over an algebraically closed field of characteristic $p>5$ is finitely generated.  In the process we prove log abundance for such pairs in the case $\kappa=2$.
\end{abstract}

\maketitle

\tableofcontents

\section{Introduction}

We will work over an algebraically closed field $k$ of characteristic $p>5$..  

The existence of good log minimal models for klt threefold pairs in characteristic $p>5$ is proven in \cite{birkar_p} when $K_X+B$ is big and generalised in \cite{birkar_waldron} to the case where only $B$ is big.  For klt threefold pairs with arbitrary boundary this would follow from \cite{birkar_p} and the log abundance conjecture.

\begin{conj}[Log abundance]\label{log-ab}
Let $(X,B)$ be a klt threefold pair over a field of characteristic $p>5$ such that $K_X+B$ is nef.  Then $K_X+B$ is semi-ample.
\end{conj}

Log abundance is proved for log canonical threefolds in characteristic zero over a sequence of papers: of Miyaoka \cite{miyaoka_chern_classes}, \cite{miyaoka_abundance} and \cite{miyaoka_kodaira_dim}, Kawamata \cite{kawamata_pluricanonical_1985} and Keel, Matsuki and McKernan \cite{keel_log_abundance}.  In positive characteristic it is still open.  A corollary of log abundance is the following, which we prove here.

\begin{thm}\label{fingen}
Let $(X,B)$ be a projective klt threefold pair over a field of characteristic $p>5$ with $\Q$-boundary $B$.  Then the log-canonical ring $$R(K_X+B)=\bigoplus H^0(\rddown{m(K_X+B)})$$ is finitely generated.
\end{thm}

Finite generation can be proved in characteristic zero by using a canonical bundle formula to reduce to lower dimension (see \cite{fujino_mori}).  This formula seems harder to obtain in positive characteristic due to inseparability and wild ramification.  Instead we prove finite generation via a special case of Conjecture \ref{log-ab}:

\begin{thm}\label{GLMM}
Let $(X,B)$ be a projective klt threefold pair over a field of characteristic $p>5$ with $\Q$-boundary $B$, such that $\kappa(X,K_X+B)=2$ and $K_X+B$ is nef.  Then $K_X+B$ is semi-ample.
\end{thm}

The proof follows that used in the proof of log abundance in \cite{flips_and_abundance}.  Using results of Keel, $K_X+B$ is endowed with map (EWM) to an algebraic space.  We run suitable LMMPs with scaling to reduce to the case where this map has equidimensional fibres.  We then show that the algebraic space is a projective variety by comparing it to a subvariety of a Chow variety, show that $\Q$-factoriality is preserved under such a morphism, and finally deduce that $K_X+B$ is the pullback of an ample divisor.

\section{Preliminaries}

\subsection{Algebraic spaces}

We need to know that some standard results about morphisms of schemes also hold for maps of algebraic spaces.  These are:

\begin{lem}[{Stein Factorisation}]\label{Stein}
Let $S$ be a scheme, and $f:X\to Y$ be a proper morphism of algebraic spaces over $S$ with $Y$ locally Noetherian.  Then $f$ can be factorised into a map with connected geometric fibres followed by a finite map.
\end{lem}
\begin{proof}
\cite[A018]{stacks}
\end{proof}

\begin{lem}[{Zariski's Main Theorem I}]\label{ZMTfin}
Let $f:X\to Y$ be a quasi-finite and separated map of algebraic spaces over a scheme $S$, such that $Y$ is quasi-compact and quasi-separated.  Then there exists a factorisation $X\to T\to Y$ so that $X\to T$ is a quasi-compact open immersion and $T\to Y$ is finite.
\end{lem}
\begin{proof}
\cite[05W7]{stacks}
\end{proof}

\begin{defn}
For a birational map of algebraic spaces $f:X\bir Y$ we define the \emph{total transform} of a subspace $Z$ is $p_2(p_1^{-1}(Z))$ where $p_1$ and $p_2$ are the projections from the closure of the graph of $f$.  We call those points in the complement of the maximal domain of definition of $f$ \emph{fundamental points}.
\end{defn}

\begin{cor}[{Zariski's Main Theorem II}]\label{ZMTbir}
Let $f:X\bir Y$ be a birational map of proper algebraic spaces.  If $P$ is a fundamental point of $f$ then its total transform $f(P)$ is connected and has positive dimension.
\end{cor}
\begin{proof}
This can be deduced from Lemma \ref{Stein} and Lemma \ref{ZMTfin} applied to the projection from the graph of the birational map.
\end{proof}

\subsection{Chow varieties}

Inseparable field extensions cause complications to arise in the theory of Chow varieties, which means that the results are weaker in positive characteristic.  This subsection summarizes the results which still hold.  The subject is treated in \cite{kollar_rational_1996}.

\begin{defn}[Algebraic cycle]
An $n$-dimensional algebraic cycle in a scheme $X$ over $k$ is a formal linear combination of $n$-dimensional reduced and irreducible subschemes of $X$ over $k$.
\end{defn}

\begin{defn}(Well defined family of algebraic cycles, \cite[I.3.10]{kollar_rational_1996})
A well defined family of algebraic cycles in $X$ consists of a reduced base scheme $Z$ and a closed subscheme $U$ of $X\times_k Z$ together with the projection morphism $g:U\to Z$ such that:
\begin{itemize}
\item{$U=\sum m_i [U_i]$ is an algebraic cycle}
\item{$g$ is proper}
\item{Every component of $U$ maps onto an irreducible component of $Z$, and every fibre is either $n$-dimensional or empty}
\item{A final technical condition must be satisfied:  we omit description of it as it is automatically satisfied when the base is normal.}
\end{itemize}
\end{defn}

There are two additional conditions on families of algebraic cycles which appear only in positive characteristic.  They are automatically satisfied in characteristic zero.

\begin{defn}\cite[I.4.7]{kollar_rational_1996}
A family of algebraic cycles satisfies the \emph{Chow-field condition} if for every $z\in Z$, the intersection of all fields of definition of the cycle corresponding to $z$ is equal to $k(z)$.

A family of algebraic cycles satisfies the \emph{field of definition condition} if for every $z\in Z$, the cycle corresponding to $z$ is defined over $k(z)$.
\end{defn}

If the base of a family is normal, that family satisfies the Chow-field condition.

These conditions lead to three attempts to define the Chow functor:

\begin{defn}\cite[I.4.11]{kollar_rational_1996}\label{chowdef}
Suppose $X$ is a scheme over a field $k$.  For every semi-normal $k$-scheme $Z$ define the following sets.
$$
Chow^{big}(X)(Z) = \left\{
\begin{aligned} &\text{Well defined proper algebraic families of non-negative}\\ &\text{cycles of } X\times Z/k\\
  \end{aligned}
\right\}
$$

$$
Chow(X)(Z) = \left\{
\begin{aligned} &\text{Well defined proper algebraic families of non-negative}\\ &\text{cycles of } X\times Z /k \text{ satisfying the Chow field condition}\\
  \end{aligned}
\right\}
$$

$$
Chow^{small}(X)(Z) = \left\{
\begin{aligned} &\text{Well defined proper algebraic families of }\\ &\text{non-negative cycles of } X\times Z /k \text{ satisfying }\\ & \text{ the field of definition condition}\\
  \end{aligned}
\right\}
$$

\end{defn}

$Chow^{\mathrm{big}}(X)(Z)$ and $Chow(X)(Z)$ agree when $Z$ is normal.  Ideally $Chow(X)$ would be a functor from the category of semi-normal $k$-schemes to the category of sets, however this is not the case.   $Chow^{small}(X)$ is such a functor, but $Chow(X)$ is only a partial functor as the defining condition is not always preserved under pullbacks. 
After a choice of very ample line bundle, the partial functors $Chow_{d,d'}^*(X)$ are defined similarly by restricting the definitions to cycles of dimension $d$ and degree $d'$.

\begin{thm}\cite[I.4.13]{kollar_rational_1996}\label{rep}
For a projective scheme $X/k$ with a given choice of very ample line bundle, there is a scheme $\Chow_{d,d'}(X)$ which coarsely represents both $Chow_{d,d'}(X)$ and $Chow_{d,d'}^{small}(X)$
\end{thm}

Although $Chow(X)$ is not a functor, it has a universal family in most cases.

\begin{thm}\cite[I.4.14]{kollar_rational_1996}\label{univ}
Let $X/k$ be a projective scheme with a given choice of very ample line bundle and $V$ an irreducible component of $\Chow_{d,d'}(X)$ of positive dimension.  Then there exists a universal family $\mathcal{U}\in Chow_{d,d'}(X)(V)$.
\end{thm}

\subsection{Reduction maps}

\begin{thm}\label{nef}
Let $X$ be a variety over an uncountable field, and $L$ a nef $\R$-divisor.  There exists a rational map $f:X\bir Z$, called the nef reduction map, with the following properties:
\begin{itemize}
\item{$f$ is proper over an open subset $U$ of $Z$}
\item{$L|_F\num 0$ on very general fibres over $U$}
\item{A curve $C$ through a very general point $x\in X$ satisfies $C\cdot L=0$ if and only if $C$ is contracted by $f$}
\end{itemize}
\end{thm}

The existence of the nef reduction map was proven for nef line bundles in characteristic zero in \cite{bauer_reduction_2002} and it was noted that the same proof applies in the more general situation above in \cite{cascini_tanaka_xu}.  When the nef reduction map exists, the dimension of its image is called the nef dimension of $L$, and is denoted $n(X,L)$.  It satisfies $\kappa(X,L)\leq\nu(X,L)\leq n(X,L)$.

\begin{defn}\cite[0.4.1]{keel_basepoint_1999}
A nef line bundle $L$ on a scheme $X$ is \emph{endowed with map (EWM)} if there is a proper map $f:X\to Z$ to an algebraic space $Z$ which contracts a subvariety $Y$ if and only if $L|_Y$ is not big.  We may always assume that such a map has geometrically connected fibres.
\end{defn}

\begin{lem}\label{ewm}
Let $X$ be a normal projective variety over an uncountable $k$ of characteristic $p>0$.  Suppose 
$L$ is a nef $\Q$-divisor on $X$ with equal Kodaira and nef dimensions $\kappa(L)=n(L)\leq 2$. 
Then $L$ is EWM to a proper algebraic space 
$V$ of dimension equal to $\kappa(L)$.
\end{lem}

\begin{proof}
\cite[7.2]{birkar_waldron}
\end{proof}

\section{Equidimensional morphisms}

\begin{defn}[Pullback of a Weil divisor]\label{pullWeil}
Let $f:X\to Y$ be a morphism of normal varieties with equidimensional fibres.  For a Weil divisor $D$ on $Y$, its pullback to $X$, denoted $D_X$, is defined as follows.  First let $Y_0$ be the smooth locus of $Y$ and $X_0$ be its pre-image by $f$.  $D|_{Y_0}$ is Cartier, so the pullback $f|_{X_0}^*(D|_{Y_0})$ is a well defined Cartier divisor on $X_0$.  The complement of $X_0$ in $X$ has codimension at least $2$, and so this pullback extends uniquely to a Weil divisor $D_X$ on $X$.
\end{defn}

\begin{lem}\label{pullback}
Let $f\colon X\to Z$ be a projective contraction between normal quasi-projective varieties over $k$
 and $L$ a nef$/Z$ $\R$-divisor on $X$ such that $L|_F\sim_\R 0$ where $F$ is the generic
fibre of $f$. Assume $\dim Z\le 3$ if $k$ has characteristic $p>0$. Then there exists a diagram
$$
\xymatrix{
X'\ar[r]^\phi\ar[d]^{f'} & X\ar[d]^f\\
Z'\ar[r]^\psi & Z
}
$$
with $\phi,\psi$ projective birational, and an $\R$-Cartier divisor $D$ on $Z'$ such that 
$\phi^* L\sim_\R f'^*D$. 
Moreover, if $Z$ is $\Q$-factorial and $f$ has equidimensional fibres, then we can take $X'=X$ and $Z'=Z$. 
\end{lem}
\begin{proof}
The first statement is proved in \cite[5.6]{birkar_waldron}, see also \cite[2.1]{kawamata_pluricanonical_1985}.
Here we prove the case where $Z$ is $\Q$-factorial and $f$ is equidimensional. 

Let $U$ be the maximal open subset of $Z$ over which $f$ is flat, and $X_U$ the base change.  $Z\backslash U$ has codimension at least $2$.  Using the argument of \cite[5.6]{birkar_waldron} (paragraph 2) and the fact that $L|_{X_U}\num 0/U$ (\cite[5.4]{birkar_waldron}), we can find a divisor $D_U$ on $U$ so that $L|_{X_U}\lin_\R f|_U^*D_U$.  Let $D$ be the closure of $D_U$ in $Z$.  Then $L\lin_\R f^*D$ because $L$ and $f^*D$ are $\R$-linearly equivalent in codimension $2$.

\end{proof}

\begin{prop}\label{QCart}
Let $f:X\to Y$ be a morphism of normal varieties over $k$ with equidimensional fibres.  Suppose $X$ is $\Q$-factorial.  Then $Y$ is also $\Q$-factorial.
\end{prop}
\begin{proof}

Let $D$ be an irreducible Weil divisor on $Y$ which we wish to show is $\Q$-Cartier.  Given a closed point $y\in Y$, let $F$ be the fibre over $y$.  We are free to replace $D$ with a multiple, so we can assume $D_X=\Div(\phi)+ A_1-A_2$ for some very ample Cartier divisors $A_1$ and $A_2$.  In particular we may assume that $A_1$ and $A_2$ do not contain any component of $F$. Then $D_X=\Div(\phi)$ on some open subvariety $U$ of $X$ containing a dense open subvariety of $F$.  Intersecting $X$ by general hyperplanes and shrinking $Y$ around $y$, we obtain a finite morphism $f|_Z:Z\to Y$ from a normal variety $Z$ such that that $D_X|_Z=D_Z=\Div(\phi)$. 

As the problem is local we may assume $Z=\Spec(A)$ and $Y=\Spec(B)$ respectively.  We have $f:\Spec(A)\to\Spec(B)$, which will be induced by a finite algebra map $\psi:B\to A$.  We may shrink further to assume $\phi$ is an element of $A$.  It then defines a $B$-linear map $\tilde{\phi}:A\to A$ by multiplication.  Define the norm $N(\phi)$ of $\phi$ to be the determinant of this linear map, which is an element of $B$.  We now claim that some multiple of $D$ is defined by $N(\phi)$.

Let $P\in\Spec(B)$, and choose $Q\in\Spec(A)$ such that $f(Q)=P$ (so that $\psi^{-1}Q=P$ as ideals).  The $B/P$ linear map $\tilde{\phi}_q:A/Q\to A/Q$ induced by multiplication by $\phi$ has determinant equal to $N(\phi) (\mathrm{mod}\ P)$.  Therefore if $N(\phi)\in P$ then $\det(\tilde{\phi})=0$ so there is some $0\neq s'\in A/Q$ such that $\tilde{\phi}_q(s')=0$.  This lifts to $s\in A\backslash Q$ such that $\phi\cdot s\in Q$.  As $Q$ is prime this implies $\phi\in Q$.  Conversely if $\phi\in Q$ then $\tilde{\phi}_q$ is the zero map and so $N(\phi)\in P$.  Thus $\Div(N(\phi))$ must be supported on $D$ and hence is equal to some multiple of $D$ as $D$ is irreducible.  

Thus we have shown that $D$ is $\Q$-Cartier in a neighbourhood of the point $y$, and $y$ was arbitrary so $D$ is $\Q$-Cartier.

\end{proof}

\section{Proof of main results}

\begin{proof}[Proof of \ref{GLMM}]

Firstly we can replace $X$ by a small crepant $\Q$-factorialisation by \cite[1.6]{birkar_p}.  Furthermore we may extend the base field to assume it is uncountable.  We include a proof of the following simple lemma for completeness.

\begin{lem}\label{n=k=2}
Let $X$ be a normal projective variety of dimension $n$ over an uncountable field with nef $\Q$-Cartier divisor $D$ such that $\kappa(X,D)= n-1$.  Then $n(X,D)=n-1$.
\end{lem}

\begin{proof}
As $\kappa=n-1$ there is a positive integer $m$ such that $mD$ is Cartier and $\phi:=\phi_{|mD|}:X\bir \PP^N$ has $n-1$ dimensional image $Y$.  Replacing $X$ by a resolution of $\phi$ and $D$ by its pullback, we may assume $\phi$ is a morphism. Sections $D'\in |mD|$ correspond to  hyperplanes in $\PP^N$ under this embedding.  Given a general point $y\in Y$, such that $\phi$ is defined on all of the fibre over $y$, let $H_1,...,H_{n-1}$ be hyperplanes intersecting at $y$.  Let $D_1,...,D_{n-1}$ be the corresponding elements of $|mD|$.  As $D$ is nef but not big, $D^n=0$.  $D_1\cdot...\cdot D_{n-1}$ is an effective $1$-cycle containing in its support the closure of the fibre over $y$.  By the nefness of $D$, $D\cdot C=0$ for each curve in this fibre.  We have covered an open subset of $X$ with $D$-trivial curves and so $n(X,D)\leq n-1$.  The lemma follows from the inequality $\kappa(X,D)\leq n(X,D)$.
\end{proof}

We may now apply Theorem \ref{nef} to obtain an almost proper nef reduction map $f:X\bir Z$ to a smooth surface $Z$.  Also by Lemma \ref{ewm}, $K_X+B$ is EWM to a proper algebraic space $V$, which we may assume is reduced and normal.  Let $g:X\to V$ be the associated map.

As $\kappa(K_X+B)=2$ we may choose $M\geq 0$ such that $K_X+B\lin_\Q M$.  $M$ is numerically trivial on a very general fibre of $f$ so it must be vertical over $Z$.  Therefore by Lemma \ref{pullback} there is a birational map $\phi:W\to X$, a contraction $f:W\to Z$ and a divisor $D$ on $Z$ such that $\phi^*(K_X+B)\lin_\Q f^*D$.  $D$ is big because $\kappa(D)=2$ and $\dim Z=2$.  We may replace $W$ and $Z$ with higher models in order to assume $W\to Z$ is flat.  Although they may now have bad singularities, replacing $D$ with its pullback still gives a well defined Cartier divisor on $Z$ satisfying the above relation. 

The following lemma is well known and we include its proof for completeness, for example see \cite[3.2]{moriwaki}.

\begin{lem}\label{contract}
There is a sequence of $K_X+B$-trivial flips and divisorial contractions leading to a model $(X',B')$ for $(X,B)$ which is also EWM to $V$, and so that no divisor on $X'$ is contracted to a point on $V$.
\end{lem}
\begin{proof}
We have the following set-up:

$$
\xymatrix{
W\ar[r]^{\phi}\ar[rd]^{f} & X\ar[r]^{g} & V \\
 & Z
}
$$

As $D$ is big we may change it up to $\Q$-linear equivalence so that $D=A+E$ where $A$ is ample and effective, $E$ is effective and the two share no components.  Also there is an effective $\Q$-divisor $M$ such that $M\lin_\Q K_X+B$ and $\phi^*M=f^*D$.

Suppose we are not already in the situation we require, so there is some reduced irreducible Weil divisor $F$ on $X$ contracted to a point by $g$, so every curve on $F$ is $K_X+B$-trivial.  Let $F_W$ be the birational transform of $F$ on $W$.  As $f$ is flat, $f(F_W)=\Gamma$ is $1$-dimensional.  By construction $\Gamma\cdot D =0$.  $A$ is ample so $\Gamma\cdot A>0$ hence $\Gamma\cdot E<0$ and so $\Gamma$ is a component of $\Supp(E)$.  Therefore $\Supp(f^*E)\subset\Supp(\phi^*M)$ contains $F_W$, so $F$ is contained in $\Supp(M)$.  

Let $\Gamma_W$ be a general curve in $F_W$ which is surjective to $\Gamma$, in particular $\Gamma_W$ should be contained in no component of $\Supp(f^*D)$ besides $F_W$.  The projection formula gives $\Gamma_W\cdot \phi^*M=0$, $\Gamma_W\cdot f^*E<0$ and $\Gamma_W\cdot f^*A>0$.  Thus there is some component $A_W$ of $\Supp(f^*A)$ such that $\Gamma_W\cap A_W\neq\emptyset$.  $A_W$  is not contracted over $X$ because every divisor which is exceptional over $X$ is contained in $\Supp(f^*E)$.  $\Supp(E)\cap\Supp(A)$ contains no divisor so as $f$ is flat, $\Supp(f^*E)\cap\Supp(f^*A)$ also contains no divisor.

The image of $\Gamma_W$ on $X$, $\Gamma_X$ is $1$-dimensional as the general curves from which we chose $\Gamma_W$ cover $F_W$ and $F_W$ is not contracted over $X$.  We know $\Gamma_X\cdot M=0$, and if $A_X$ is the birational transform of $A_W$ then $\Gamma_X\cdot A_X>0$.  Therefore we must have $\Gamma_X\cdot F<0$.

Run a $K_X+B+\epsilon F$-MMP with scaling of some ample divisor for $\epsilon$ sufficiently small, which terminates by \cite[1.6]{birkar_waldron}.  We claim that every step is $K_X+B$-trivial.  This holds for the first step because any curve with $K_X+B+\epsilon F<0$ is contained in $F$ and so is contracted over $V$.  This implies that $K_X+B$ is the pullback of a $\Q$-Cartier divisor on the contracted variety, and hence so is the new $K_{X'}+B'$ after the contraction or flip.  Hence $K_{X'}+B'$ is also nef, and is still EWM to $V$.  The birational transform of $F$ is still contracted over $V$.  Thus the same argument applies inductively to each step.

The LMMP terminates on some model $(X',B')$, and $K_{X'}+B'$ is EWM to V.  The divisor $F$ must have been contracted during the LMMP, as otherwise its birational transform $F'$ would be covered by $K_{X'}+B'+\epsilon F'$-negative curves.

After repeating this procedure finitely many times we arrive at the model described in the statement.
\end{proof}

Replace $(X,B)$ with the pair $(X',B')$ constructed in Lemma \ref{contract} to assume $g:X\to V$ has equidimensional fibres.  Next we now prove that $V$ is projective by comparing it with a subvariety of a Chow variety.  For an example of Chow varieties in a similar context, see \cite{flips_and_abundance}.

\begin{lem}
$V$ is a projective variety.
\end{lem}

\begin{proof} \

\textit{Step 1}: We construct a projective variety $Z^\nu$ birational to $V$ from a subvariety of a Chow variety. \smallskip

By \cite[0ADD]{stacks} an algebraic space is a scheme in codimension $1$.  Therefore there is a smooth open subvariety $U$ of $V$ such that $V\backslash U$ has codimension at least $2$, for which $g_U:X_U=g^{-1}(U)\to U$ is a morphism of varieties.  As $V$ is reduced, $g_U$ is a well defined family of algebraic cycles.  It is therefore an element of $Chow^{\mathrm{big}}(U)$.  As $U$ is normal it is also an element of $Chow(U)$, which is enough to ensure that there is a morphism $U\to \Chow(X)$ by Theorem \ref{rep}.  $U$ is separated and $\Chow(X)$ is proper, therefore $U\to\Chow(X)$ is separated.  

As each point of $U$ represents a cycle supported on a different locus of $X$, the image of $U$ is contained in a component $C$ of $\Chow(X)$ of positive dimension.  By Theorem \ref{univ} there is a universal family $h:\mathcal{U}\to C$.  Let $Z$ be the closure of the image of $U$ in $\Chow(X)$ with reduced scheme structure, and $X'$ the inverse image of $Z$ in $\mathcal{U}$.  The morphism $X'\to Z$ has geometrically connected fibres.  The morphism $U\to Z$ is quasi-finite as every fibre over a $k$ point contains either $1$ or $0$ closed points (the cycle represented by a point on $Z$ is supported on the fibre of at most one point in $U$).  Apply Zariski's Main Theorem \ref{ZMTfin} to this morphism, to factorise it as $U\to \hat{Z}\to Z$ where $U\to \hat{Z}$ is an open immersion and $\hat{Z}\to Z$ is finite.  Let $\hat{X}=X\times_z\hat{Z}$.  Finally let $Z^\nu$ be the normalisation of $\hat{Z}$ and $X^\nu = Z^\nu \times_{\hat{Z}}\hat{X}$.  We get an open immersion $U\to Z^\nu$ and $X^\nu\to Z^\nu$ with $1$-dimensional fibres and purely $3$-dimensional $X^\nu$.  The universal family $\mathcal{U}$ is a subvariety of $X\times_k C$, so comes with a natural projection to $X$.  By composition we get a natural projection $\pi:X^\nu\to X$.

$$
\xymatrix{
X &\mathcal{U}\ar[l]^\pi\ar[d]^h & X'\ar[l]\ar[d]^{h'} & X^\nu\ar[l]\ar[d]^{h^\nu}& X_U \ar[l]\ar[r]\ar[d] & X\ar[d]^{g}  \\
& C & Z\ar[l] & Z^\nu \ar[l]^\nu& U \ar[l]\ar[r] & V
}
$$

$Z^\nu$ and $V$ are birational as they both contain $U$ as an open subspace.\medskip

\textit{Step 2:} We apply Zariski's Main Theorem (\ref{ZMTbir}) to show that $V\isom Z^\nu$.  \smallskip

As $V\backslash U$ is of codimension at most $2$, $Z^\nu\bir V$ has no fundamental points as the total transform of such a point would be of codimension $1$ in $V\backslash U$.  So it is enough to show that there are no fundamental points of $V\bir Z^\nu$.

Suppose for contradiction that $P\in V$ is a fundamental point of $V\bir Z^\nu$ with $1$-dimensional total transform $T\subset Z^\nu$.  Let $\Gamma$ be a general curve through $P$.  We may assume $\Gamma$ passes through no other point of $V\backslash U$.  Let $\Gamma^\nu$ be the closure of $\Gamma|_U$ in $Z^\nu$.  As $\Gamma$ is general, the intersection of $\Gamma^\nu$ with $T$ consists of a general finite set of points of $T$.  

We now claim that one of these points has image in $C$ representing a cycle with support contained in $g^{-1}(P)$.  
Pull back the universal family from $X^\nu\to Z^\nu $ to $Y^\nu\to \Gamma^\nu$.  $Y^\nu$ is a geometrically connected $2$-dimensional cycle, and all but a $1$-dimensional subvariety is over $U$.  Let $Y\subset X$ be the pre-image of $\Gamma$ by $g$.

The image of $Y^\nu$ on $X$ via the universal family is supported on the closure of the image of $X_U\cap Y^\nu$, and this is equal to $\overline{X_U}=g^{-1}(\Gamma)=Y$.  The image of $X_U\cap Y^\nu$ is supported precisely  on the union of the supports of the cycles corresponding to the points of $\Gamma_U$.  The supports of these cycles cover all of $Y$ except for $g^{-1}(P)$.  We know that the image of the projection from $Y^\nu$ to $X$ must be supported on all of $Y$.  Therefore any component of the support of $g^{-1}(P)$ must be contained in the support of some cycle represented by a point $Q\in \Gamma^\nu\backslash \Gamma_U$.  Suppose that the cycle represented by $Q$ also contains some other $1$-dimensional subscheme of $Y$ not contained in $g^{-1}(P)$.  This is impossible as that subscheme could not be contained in any $2$-dimensional component of $Y^\nu$ (as each of the cycles near it are disjoint from one another but supposedly intersect this cycle).   Thus we have shown that $Q$ represents a cycle supported only on the support of $g^{-1}(P)$. 

We have proved the claim that when a general curve intersects $T$ in a finite set of points, one of these points represents a cycle supported in $g^{-1}(P)$.  But there are only countably many $1$-dimensional cycles supported there, and these can be represented by only countably many points of $C$.  The pre-images of these in $Z^\nu$ are also countably many points.  We can find some general curve $\Gamma^\nu$ intersecting $T$ away from these points.  This is a contradiction to the existence of the point $P$ and so $V\isom Z^\nu$.

\end{proof}

$X$ is $\Q$-factorial and $g$ has equidimensional fibres so $V$ is $\Q$-factorial by Proposition \ref{QCart}.  Also, by Lemma \ref{pullback},  $K_X+B\lin g^*D$ for some big and nef divisor $D$.  $K_X+B$ is EWM with map $g$, so $D$ is big and strictly nef on a surface.  Therefore Kodaira's lemma and the Nakai-Moishezon criterion imply that $D$ is ample.
\end{proof}

\begin{proof}[Proof of \ref{fingen}]
The case $\kappa=3$ is proved in \cite[1.3]{birkar_p}, $\kappa=2$ follows from Theorem \ref{GLMM}, $\kappa=1$ from Lemma \ref{k=1} below, and $\kappa=0$ and $\kappa = -\infty$ are obvious.

\begin{lem}\label{k=1}
Let $X$ be a projective normal variety and $L$ a $\Q$-Cartier $\Q$-divisor on $X$ with $\kappa(X,L)=1$.  Then the section ring $\oplus_k H^0(X,\rddown{kL})$ is finitely generated.
\end{lem}
\begin{proof}
It is enough to show that $R = \oplus_kH^0(X,\rddown{kaL})$ is finitely generated for any $a$, so we may freely replace $L$ with a higher power.  By doing so we may assume $L$ is an effective Cartier divisor and that the following natural inclusions are strict. 
$$H^0(X,\OO_X)\subset H^0(X,L)\subset H^0(X,2L)\subset...$$ 

$\mathcal{B}^0_0=\{1\}$ is a $k$ basis for $H^0(X,\OO_X)$, where $1$ is the ring identity element contained in $H^0(X,\OO_X)$.  Let $1_L$ be the image of $1$ in $H^0(X,L)$ under this inclusion.  Let $\mathcal{B}^0_1=\{1_L\}\subset H^0(X,L)$, and let $\mathcal{B}^1_1$ be a collection of functions completing $\mathcal{B}^0_1$ to a basis of $H^0(X,L)$.  Fix some $x\in \mathcal{B}^1_1$.  Inductively we let $\mathcal{B}^j_{i+1}$ be the vectors in $1_L\cdot \mathcal{B}^j_i$ for each $j\in\{0,...,i\}$, and take $\mathcal{B}^{i+1}_{i+1}$ to complete $\cup_{j=0}^i\mathcal{B}^j_{i+1}$ to a basis of $H^0(X,(i+1)L)$.  We claim that we may choose $\mathcal{B}^{i+1}_{i+1}$ to contain $x\cdot\mathcal{B}_i^i$.  This is equivalent to saying that $$\mathrm{Span}(\{ x\cdot\mathcal{B}_i^i\})\cap\mathrm{Span}( \{{1_L}\cdot\mathcal{B}_i\})=0$$ and $x\cdot \mathcal{B}^i_i$ is linearly independent.  For the first assertion, there can be no such non-zero element by comparing the order of the zero at $L$ implied by membership of the two subspaces.  The second assertion follows from the linear independence of $\mathcal{B}^i_i$ and that $K(X)$ is a domain.

By construction, multiplication by $x\in\mathcal{B}^1_1$ induces an injection $\mathcal{B}^k_k\to\mathcal{B}_{k+1}^{k+1}$.  This implies that $|\mathcal{B}_k|-|\mathcal{B}_{k-1}|\leq |\mathcal{B}_{k+1}|-|\mathcal{B}_k|$, and so the sequence of integers $|\mathcal{B}_k|-|\mathcal{B}_{k-1}|$ is increasing.

The condition $\kappa(L)=1$ implies that by possibly replacing $L$ with a higher multiple there are constants $A$ and $B$ such that $$Ak\leq H^0(X,kL)=\sum_{j=0}^k(|\mathcal{B}_j|-|\mathcal{B}_{j-1}|)+1\leq Bk$$
This implies that the sequence stabilises at an integer at most $B$, and so the injections $x:\mathcal{B}_k^k\to\mathcal{B}_{k+1}^{k+1}$ are bijections for $k$ large enough.  

\end{proof}

\end{proof}

\section*{Acknowledgements}
The author would like to thank Caucher Birkar for his advice and support and the referee for their helpful comments.  This research was funded by an EPSRC studentship.

\bibliography{bib}{}
\bibliographystyle{abbrv}
\end{document}